\documentclass[12pt]{amsart}

\usepackage{amsmath,amssymb,epic,lscape}

\newtheorem{theorem}{Theorem}[section]
\newtheorem{proposition}[theorem]{Proposition}

\theoremstyle{definition}

\theoremstyle{remark}
\newtheorem{remark}[theorem]{Remark}

\numberwithin{equation}{section}

\def\DJ{{\hbox{{\thinspace}D\kern-.8em\raise.15ex\hbox{--}\kern.35em}}}
\def\DJo{{$\;$\kern-.42em{\DJ}okovi\'c}}

\renewcommand{\subjclassname}{\textup{2000} Mathematics Subject
Classification}

\font\germ=eufm10

\def\al{{\alpha}}
\def\be{{\beta}}
\def\ga{{\gamma}}
\def\vf{{\varphi}}

\def\bR{{\mbox{\bf R}}}
\def\bZ{{\mbox{\bf Z}}}
\def\bC{{\mbox{\bf C}}}

\def\bH{{\mbox{\bf H}}}
\def\bx{{\mbox{\bf x}}}
\def\by{{\mbox{\bf y}}}
\def\pA{{\mathcal A}}

\def\pP{{\mathcal P}}

\def\Cb{{\bar{C}}}

\def\diag{{\rm diag}}
\def\tr{{\rm tr\;}}

\def\GL{{\mbox{\rm GL}}}
\def\SL{{\mbox{\rm SL}}}
\def\SO{{\mbox{\rm SO}}}
\def\Sp{{\mbox{\rm Sp}}}

\def\Ort{{\mbox{\rm O}}}

\def\gg{{\mbox{\germ g}}}
\def\gl{{\mbox{\germ gl}}}
\def\gk{{\mbox{\germ k}}}
\def\gp{{\mbox{\germ p}}}
\def\sp{{\mbox{\germ sp}}}

\begin{document}

\title[Symplectic invariants of matrices]
{Symplectic polynomial invariants of one or two matrices
of small size}

\author[D.\v{Z}. \DJ okovi\'{c}]
{Dragomir \v{Z}. \DJ okovi\'{c}}

\address{Department of Pure Mathematics, University of Waterloo,
Waterloo, Ontario, N2L 3G1, Canada}

\email{djokovic@uwaterloo.ca}

\thanks{
The author was supported by an NSERC Discovery Grant.}


\date{}

\begin{abstract}
The algebra of holomorphic polynomial $\Sp_{2n}$-invariants of $k$ 
complex $2n\times 2n$ matrices (under diagonal conjugation action) 
is generated by the traces of words in these matrices and their 
symplectic adjoints. No concrete minimal generating set is known 
for this algebra apart from the cases $n=1$, when $\Sp_2=\SL_2$, 
and $n=2$, $k=1$. We construct such sets in the cases $n=k=2$ and 
$n=3$, $k=1$. In the latter case we also construct a homogeneous system 
of parameters and a Hironaka decomposition of the algebra.
\end{abstract}

\maketitle
\subjclassname{ 13A50, 14L35, 20G20 }

\section{Introduction} \label{Uvod}

We shall assume throughout this paper that the underlying field 
is $\bC$, the field of complex numbers.
For the basic results of the invariant theory of classical groups
the reader may wish to consult \cite{GW,HK,SM,PO,HW} or the very
recent book \cite[Chapter 11]{CP2}.
While the First and Second Fundamental Theorems (FFT and SFT)
for invariants of several, say $k$, matrices under the simultaneous 
conjugation action by one of the classical groups
($\SL_n$, $\Ort_n$, $\SO_n$ or $\Sp_{2n}$) are known,
the concrete results concerning the algebra of polynomial invariants
such as the Poincar\'{e} series,
a minimal set of homogeneous generators (MSG),
a homogeneous system of parameters (HSOP),
the Hironaka decomposition, ...
is known only in a handful of low-dimensional cases.
We recall that Schwarz \cite{GS} has determined all representations 
of simple Lie groups having a regular algebra of invariants, 
i.e., the algebra of invariants isomorphic to a polynomial algebra.
Moreover in each of these cases he has constructed an MSG.

This paper deals exclusively with the symplectic group $\Sp_{2n}$ 
(see the next section for its matrix definition).
It is well known that the algebra of $\Sp_{2n}$-invariants of $k$ 
complex $2n\times 2n$ matrices, under diagonal conjugation action
\begin{equation} \label{dij-akcija}
a\cdot(x_1,\ldots,x_k)=(ax_1a^{-1},\ldots,ax_ka^{-1}),
\end{equation}
is generated by the traces of words in the matrices $x_i$ and 
their symplectic adjoints $x_i^*$.
No concrete MSG is known for any $k\ge2$ apart from the case $n=1$ 
when $\Sp_2=\SL_2$. We shall construct such a set consisting of 
136 tracesin the case $n=k=2$. In the case $n=3$, $k=1$ we 
construct an HSOP consisting of 15 traces and an MSG consisting 
of 28 traces.

Let $\Sp_{2n}\subseteq\GL_{2n}$ be the standard Lie group embedding 
and let $\gk=\sp_{2n}\subseteq\gg=\gl_{2n}$ be the induced embedding 
of their Lie algebras. We can view $\gg$ as an $\Sp_{2n}$-module via 
the restriction of the adjoint representation of $\GL_{2n}$.
There is a unique direct decomposition $\gg=\gk\oplus\gp$ of this 
module. 

More generally, we consider the direct sum 
$$
V=k_1\gp \oplus k_2\gk
$$
of $k_1$ copies of $\gp$ and $k_2$ copies of $\gk$. 
In the case $k_1=k_2=k$ this is just the direct sum 
of $k$ copies of the complex matrix algebra $M_{2n}$ with the diagonal 
action (\ref{dij-akcija}) of $\Sp_{2n}$. Let $\pP=\pP_V$ denote the 
algebra of complex holomorphic polynomial functions on $V$.
By fixing the above direct decomposition of $V$, we can introduce a (non-canonical) $\bZ^{k_1+k_2}$-gradation on $\pP$. For that 
purpose let us denote the standard generators of $\bZ^{k_1+k_2}$ 
by $e_r$, $1\le r\le k_1+k_2$.
Then we assign to the coordinate functions of the $r$th copy of 
$\gp$ the degree $e_r$ and to those of the $s$th copy of $\gk$ the
degree $e_{k_1+s}$. For simplicity, 
we shall refer to this gradation simply as the \emph{multigradation}.
Let $\Cb$ denote the subalgebra of $\Sp_{2n}$-invariant
functions in $\pP$. It is easy to see that $\Cb$ inherits the
multigradation from $\pP$.

Our first objective is to compute the Poincar\'{e} series of the 
algebra $\Cb$ for some small values of $k_1$ and $k_2$. 
This has been done by Berele and Adin \cite{BA} when $n=2$ and 
$k_1+k_2=2$. In section \ref{n=2} we extend their computations to 
cover the cases $k_1+k_2=3$. In all cases we used MAPLE \cite{Maple} 
and the well known Molien--Weyl formula (see \cite{DK}) to compute 
the Poincar\'{e} series.

Our second objective is to construct an MSG for some of the 
algebras $\Cb$. In section \ref{2-mat-4} we consider the case 
$n=2$ and $k_1=k_2=2$, i.e., the case of symplectic invariants of 
two complex $4\times4$ matrices $x$ and $y$. While we were not 
able to compute the multigraded Poincar\'{e} series of $\Cb$, 
we have done that for a suitable $\bZ^2$-gradation. 
We have also constructed an MSG of $\Cb$.
It consists of 136 traces of words in $x$, $y$ and their 
symplectic adjoints $x^*$, $y^*$. The maximum (total) degree 
of these generators is 9.

In section \ref{n=3} we give the bigraded Poincar\'{e} series of 
$\Cb$ in the cases $n=3$ and $k_1+k_2=2$, as well as for $k_1=3$, 
$k_2=0$.

In section \ref{1-mat-6} we analyze the algebra $\Cb$ when $n=3$ and 
$k_1=k_2=1$, i.e., the algebra of symplectic invariants of a single
complex matrix $z\in M_6$. We use the $\bZ$-gradation and construct 
an MSG consisting of 28 traces, $\tr w(z,z^*)$. 
We also propose a candidate for an HSOP in this setting.

In section \ref{Hir} we investigate the same algebra $\Cb$ as in
section \ref{1-mat-6} but now we write $z=x+y$ with $x^*=x$ and 
$y^*=-y$ and express the elements of $\Cb$ as polynomials in the 
traces $\tr w(x,y)$. This setting provides $\Cb$ with a 
$\bZ^2$-gradation. By using this bigradation, we construct an HSOP 
of $\Cb$ (of size 15) and a new MSG (of size 28). Next we construct 
a Hironaka's decomposition of $\Cb$ with respect to the polynomial 
subalgebra $\Cb^\#$ generated by our HSOP, i.e., we find a basis of 
$\Cb$ considered as a $\Cb^\#$-module (which is free of rank 36). 
Finally, we compute the unique syzygy of our MSG of minimal degree 
(degree 14). Its bidegree is (6,8).

For the known results of this kind for the classical groups 
$\SL_n$ and $\SO_n$ the reader may consult our recent papers 
\cite{DZ1,DZ2} and the references given there.

\section{Preliminaries} \label{Objasnjenja}

Let $M_{k}=M_k(\bC)$ denote the algebra of complex $k\times k$ 
matrices. Let $X^T$ denote the transpose of $X\in M_k$,
and let $I_k\in M_k$ denote the identity matrix. 
Denote by $\theta$ the involutory automorphism of the general
linear group $G=\GL_{2n}(\bC)$ whose fixed point set is the
symplectic group $\Sp_{2n}$. The differential ${\rm d}\theta$
of $\theta$ at $I_{2n}$ is an involutory automorphism
of the Lie algebra $\gg=\gl_{2n}=M_{2n}$ of $G$. The 
fixed point set of ${\rm d}\theta$ on $\gg$ is the Lie algebra
$\gk=\sp_{2n}$ of $\Sp_{2n}$. The other eigenspace of $\theta$,
for eigenvalue $-1$, will be denoted by $\gp$.
For $X\in\gg$ set $X^*=-{\rm d}\theta(X)$. We refer to $X^*$ as the
{\em symplectic adjoint} of $X$. The map $X\to X^*$
is an involutory anti-automorphism (of symplectic type) of the
associative algebra $M_{2n}$. Then $\gp$ resp. $\gk$ is the
set of symmetric resp. skew-symmetric elements of $M_{2n}$
with respect to this involution.

To be concrete, let us define the symplectic group by
$$ \Sp_{2n}=\{A\in\GL_{2n}:A^TJA=J\},\qquad
J=\left[ \begin{array}{rr}0&-I_n\\I_n&0\end{array} \right]. $$
Then we have
\begin{eqnarray*}
&& \theta(X)=J^{-1}(X^T)^{-1}J,\quad X\in\GL_{2n}; \\
&& {\rm d}\theta(X)=-J^{-1}X^TJ,\quad X\in M_{2n}, \\
&&  X^*=J^{-1}X^TJ,\quad X\in M_{2n}, \\
&& \gk=\{X\in M_{2n}:X^*=-X\}, \\
&& \gp=\{X\in M_{2n}:X^*=X\}.
\end{eqnarray*}
If we partition $X$ into four $n\times n$ blocks, then
$$ X=\left[\begin{array}{rr}X_1&X_2\\X_3&X_4\end{array}\right],\quad
X^*=\left[\begin{array}{rr}X_4^T&-X_2^T\\-X_3^T&X_1^T\end{array} 
\right]. $$
Hence
\begin{eqnarray} \label{def-k}
\gk &=& \left\{ \left[ \begin{array}{rr} X_1 & X_2 \\ X_3 & -X_1^T 
\end{array} \right]~:~X_2^T=X_2,~X_3^T=X_3 \right\}, \\
\label{def-p}
\gp &=& \left\{ \left[ \begin{array}{rr} X_1 & X_2 \\ X_3 & X_1^T 
\end{array}\right]~:~X_2^T=-X_2,~X_3^T=-X_3 \right\}.
\end{eqnarray}
Note that $\dim \gk=n(2n+1)$ and $\dim \gp=n(2n-1)$.

As $\Sp_{2n}$-module, $\gk$ is simple (the adjoint module) while 
$\gp$ is the direct sum of the 1-dimensional trivial module 
spanned by $I_{2n}$ and another simple module, 
the second fundamental module of $\Sp_{2n}$.

Recall that we use the direct decomposition $V=k_1\gp\oplus k_2\gk$ 
to define the multigradation of $\pP$ and $\Cb$.
For brevity, let us write an element of $\bZ^{k_1+k_2}$ as
$(\al,\be)$, where $\al=(\al_1,\ldots,\al_{k_1})$ and
$\be=(\be_1,\ldots,\be_{k_2})$. Let 
$x_1,\ldots,x_{k_1}$ and $y_1,\ldots,y_{k_2}$ 
be independent commuting variables. We set
$\bx=(x_1,\ldots,x_{k_1})$ and $\by=(y_1,\ldots,y_{k_2})$,
and define
$$ \bx^\al=x_1^{\al_1}\cdots x_{k_1}^{\al_{k_1}},
\quad \by^\be=y_1^{\be_1}\cdots y_{k_2}^{\be_{k_2}}. $$
By definition, the multigraded Poincar\'{e} series of
$\Cb$ is given by:
$$ P(\Cb;\bx,\by)=\sum_{(\al,\be)\ge0} d_{\al,\be}\bx^\al \by^\be $$
where $d_{\al,\be}$ is the dimension of the homogeneous component
of $\Cb$ of multidegree $(\al,\be)$. 
The case when $k_1=k_2=k$ is of special interest. In that case we
refer to $\Cb$ as the {\em algebra of symplectic invariants of $k$ 
matrices}.

The problem of computing these multigraded Poincar\'{e} series
has been treated by A. Berele and R. Adin in \cite[Section 5]{BA}.
(We have borrowed some of the notation from that paper.)
Let us recall some of their results.

First of all they have shown that the above Poincar\'{e} series is
the Taylor expansion at the origin of a rational function in the
variables $\bx$ and $\by$. We shall write this rational function as
$$ P(\Cb;\bx,\by)=\frac{N(\Cb;\bx,\by)}{D(\Cb;\bx,\by)}
=\frac{N(k_1,k_2)}{D(k_1,k_2)}, $$
where $N$ and $D$ are polynomials normalized so that they take
value 1 at the origin. Moreover we require $N$ and $D$ to be
relatively prime. Then these numerators and denominators are
uniquely determined.

Second, they have shown (see their Theorem 9) that the Poincar\'{e}
series satisfies the functional equation:
$$ P(\Cb;\bx^{-1},\by^{-1})=(-1)^{d}(x_1\cdots x_{k_1})^{n(2n-1)}
(y_1\cdots y_{k_2})^{n(2n+1)}P(\Cb;\bx,\by), $$
where $d=n(k_1+k_2+1)$ and 
$$ \bx^{-1}=(x_1^{-1},\ldots,x_{k_1}^{-1}),\quad
\by^{-1}=(y_1^{-1},\ldots,y_{k_2}^{-1}). $$

Third, in their Lemma 25 they give the multigraded Poincar\'{e}
series for $n=2$ and $k_1+k_2=2$. We note that there is a misprint
in the case $k_1=0$, $k_2=2$: The term $+y_1^4 y_2^4$ is missing
from the numerator.

Procesi \cite{CP1} has shown that in the case $k_1=k_2=k$ the algebra
$\Cb$ of $\Sp_{2n}$-invariants of $k$ generic $2n\times 2n$ matrices
is generated by the traces of words in the $k$ generic matrices
and their symplectic adjoints of length at most $2^{2n}-1$.
(In fact his bound was written incorrectly as $2^n-1$ but it is clear 
from the context that the intended bound is $2^{2n}-1$.)
A better bound, $(2n)^2$, follows from a theorem of Razmyslov
\cite{YR}. The problem of finding the best possible bound is closely
related to Kuzmin's conjecture \cite{EK,DF} about associative nil algebras
of given exponent, say $m$. Namely, he conjectured that such algebra
is nilpotent of exponent $m(m+1)/2$. If true, this conjecture in the
case $m=2n$ would give the bound $n(2n+1)$. The conjecture has
been proved for $m\le4$, see \cite{MVL} for the case $m=4$.

By $\bH$ we denote the division algebra of real quaternions and
by $M_n(\bH)$ the real algebra of $n\times n$ matrices over $\bH$.
Occasionally we use the standard quaternionic units $1,i,j,k$,
which will be evident from the context. The trace of a quaternion
$q$ is defined by $\tr q=q+\bar{q}=2\Re(q)$. 
For $z=(z_{ij})\in M_n(\bH)$, the trace is defined by
$$ \tr z=\tr\left( \sum_{i=1}^n z_{ii} \right). $$
The adjoint, $z^*$, of $z$ is defined as the conjugate
transpose of $z$.

\section{Some Poincar\'{e} series $P(\Cb;\bx,\by)$ for $n=2$} \label{n=2}

For $n=2$ and $k_1=k_2=1$, $\Cb$ is the algebra of
symplectic invariants of a single $4\times 4$ complex matrix $x$.
This algebra is in fact a polynomial algebra in six generators,
the traces of the six matrices:
$$ x,\; x^2,\; x^3,\; x^4,\; xx^*,\; x^*x^3. $$
For this result see \cite{GS} or \cite[Corollary 6.1]{AB}. 
One can easily show that
the matrix $x^*x^3$ can be replaced by $(xx^*)^2$ or $x^2(x^*)^2$.
We point out that \cite[Theorem 6.1]{AB} has a misprint:
In the formula for the Poincar\'{e} series, the factor 
$(1-t^3)^{-2}$ should be replaced with $(1-t^3)^{-1}$.

We present in this section the results of our computations of the
multigraded or simply graded Poincar\'{e} series for the algebra
$\Cb$ of $\Sp_4$-invariants of the direct sum $V$ of $k_1$ copies of 
$\gp$ and $k_2$ copies of $\gk$ when $k_1+k_2=3$.

Since $B_2=C_2$ in Lie theory, it is evident that $\gp$ can be viewed
as $\SO_5$-module. As such it is a direct sum of a 1-dimensional
trivial module and the defining 5-dimensional module. 
Hence the FFT and SFT for the direct sums of
copies of $\gp$ is known from classical invariant theory for
the special orthogonal groups \cite{CP2}. By using \cite{GS},
we see that, in the case $k_2=0$, $\Cb$ is a polynomial algebra
iff $k_1\le5$. In these cases it is easy to write down an MSG
for $\Cb$.

We have mentioned that the cases $k_1+k_2=2$ have been dealt with in \cite{BA}.
When $k_1=0$, $k_2=2$ the $\bZ$-graded Poincar\'{e} series has been computed
previously in \cite[Section 8.2]{BB}.

Assume now that $k_1+k_2=3$.
The formulae for $N$ and $D$ are simple for $k_2\le1$:

\begin{eqnarray*}
N(3,0) &=& 1, \\
D(3,0) &=& (1-x_1)(1-x_2)(1-x_3)(1-x_1^2)(1-x_2^2)(1-x_3^2)\cdot \\
&& (1-x_1x_2)(1-x_2x_3)(1-x_1x_3), \\
N(2,1) &=& (1+x_1x_2y_1^2)(1+x_1x_2y_1^3), \\
D(2,1) &=& (1-x_1)(1-x_2)(1-x_1^2)(1-x_1x_2)(1-x_2^2)(1-y_1^2)\cdot \\
&& (1-x_1x_2y_1)(1-x_1y_1^2)(1-x_2y_1^2)\cdot \\
&& (1-x_1^2y_1^2)(1-x_2^2y_1^2)(1-y_1^4).
\end{eqnarray*}

The numerator $N(1,2)$ can be written as
\begin{eqnarray*}
N(1,2) &=& c_0+c_1y_1y_2x_1+c_2y_1y_2x_1^2+c_3y_1^2y_2^2(y_1+y_2)x_1^3 \\
&& +c_4y_1^3y_2^3x_1^4+c_5y_1^4y_2^4x_1^5+c_6y_1^6y_2^6x_1^6,
\end{eqnarray*}
where
\begin{eqnarray*}
c_0 &=& 1-y_1^{2}y_2-y_1y_2^{2}+y_1^{2}y_2^{2}+y_1^{4}y_2^{2}
+2y_1^{3}y_2^{3}+y_1^{2}y_2^{4}  \\
&& +y_1^{4}y_2^{4}-y_1^{5}y_2^{4}-y_1^{4}y_2^{5}+y_1^{6}y_2^{6} \\
c_1 &=& y_1+y_2+y_1^{2}+2y_1y_2+y_2^{2}
+y_1^{3}+2y_1^{2}y_2+2y_1y_2^{2}+y_2^{3} \\
&& -y_1^{2}y_2^{2}-y_1^{5}y_2-2y_1^{4}y_2^{2}
-2y_1^{3}y_2^{3}-2y_1^{2}y_2^{4}-y_1y_2^{5} \\
&& +y_1^{4}y_2^{3}+y_1^{3}y_2^{4}+y_1^{6}y_2^{3} 
+2y_1^{5}y_2^{4}+2y_1^{4}y_2^{5}+y_2^{6}y_1^{3} \\
&& -y_1^{7}y_2^{5}-y_1^{6}y_2^{6}-y_1^{5}y_2^{7} \\
c_2 &=& y_1+y_2+y_1^{2}+2y_1y_2+y_2^{2} 
+y_1^{3}+2y_1^{2}y_2+2y_1y_2^{2}+y_2^{3} \\
&& -y_1^{2}y_2^{2}-y_1^{5}y_2-y_1^{4}y_2^{2}-y_1^{3}y_2^{3}
-y_1^{2}y_2^{4}-y_1y_2^{5}+y_1^{4}y_2^{3}  \\
&& +y_1^{3}y_2^{4}-y_1^{7}y_2^{3}-y_1^{6}y_2^{4}-2y_1^{5}y_2^{5}
-y_1^{4}y_2^{6}-y_1^{3}y_2^{7}-y_1^{7}y_2^{4} \\
&& -2y_1^{6}y_2^{5}-2y_1^{5}y_2^{6}-y_1^{4}y_2^{7}
-y_1^{7}y_2^{5}-y_1^{6}y_2^{6}-y_1^{5}y_2^{7} \\
&& +y_1^{8}y_2^{6}+y_1^{7}y_2^{7}+y_1^{6}y_2^{8} \\
c_3 &=& y_1y_2+y_1^{3}+2y_1^{2}y_2+2y_1y_2^{2}+y_2^{3}
+y_1^{3}y_2+y_1y_2^{3}+y_1^{4}y_2 \\
&& +y_1y_2^{4}-y_1^{5}y_2-3y_1^{4}y_2^{2}
-3y_1^{3}y_2^{3}-3y_1^{2}y_2^{4}-y_1y_2^{5} \\
&& -2y_1^{5}y_2^{2}-2y_1^{4}y_2^{3}-2y_1^{3}y_2^{4}
-2y_1^{2}y_2^{5}-y_1^{6}y_2^{2}-3y_1^{5}y_2^{3} \\
&& -3y_1^{4}y_2^{4}-3y_1^{3}y_2^{5}-y_1^{2}y_2^{6}
+y_1^{6}y_2^{3}+y_1^{3}y_2^{6}+y_1^{6}y_2^{4} \\
&& +y_1^{4}y_2^{6}+y_1^{7}y_2^{4}+2y_1^{6}y_2^{5}
+2y_1^{5}y_2^{6}+y_1^{4}y_2^{7}+y_1^{6}y_2^{6}, \\
c_4 &=& y_1^{2}+y_1y_2+y_2^{2}-y_1^{3}y_2-y_1^{2}y_2^{2}
-y_1y_2^{3}-y_1^{4}y_2-2y_1^{3}y_2^{2} \\
&& -2y_1^{2}y_2^{3}-y_1y_2^{4}-y_1^{5}y_2-y_1^{4}y_2^{2}
-2y_1^{3}y_2^{3}-y_1^{2}y_2^{4}-y_1y_2^{5} \\
&& +y_1^{5}y_2^{4}+y_1^{4}y_2^{5}-y_1^{7}y_2^{3}
-y_1^{6}y_2^{4}-y_1^{5}y_2^{5}-y_1^{4}y_2^{6} \\
&& -y_1^{3}y_2^{7}-y_1^{6}y_2^{6}+y_1^{8}y_2^{5}
+2y_1^{7}y_2^{6}+2y_1^{6}y_2^{7}+y_1^{5}y_2^{8} \\
&& +y_1^{8}y_2^{6}+2y_1^{7}y_2^{7}+y_1^{6}y_2^{8}
+y_1^{8}y_2^{7}+y_1^{7}y_2^{8}, \\
c_5 &=& -y_1^{2}-y_1y_2-y_2^{2}+y_1^{4}y_2+2y_1^{3}y_2^{2}
+2y_1^{2}y_2^{3}+y_1y_2^{4} \\
&& +y_1^{4}y_2^{3}+y_1^{3}y_2^{4}-y_1^{6}y_2^{2}
-2y_1^{5}y_2^{3}-2y_1^{4}y_2^{4}-2y_1^{3}y_2^{5} \\
&& -y_1^{2}y_2^{6}-y_1^{5}y_2^{5}+y_1^{7}y_2^{4}
+2y_1^{6}y_2^{5}+2y_1^{5}y_2^{6}+y_1^{4}y_2^{7} \\
&&  +y_1^{7}y_2^{5}+2y_1^{6}y_2^{6}
+y_1^{5}y_2^{7}+y_1^{7}y_2^{6}+y_1^{6}y_2^{7}, \\
c_6 &=& 1-y_1^{2}y_2-y_1y_2^{2}+y_1^{2}y_2^{2}+y_1^{4}y_2^{2}
+2y_1^{3}y_2^{3}+y_1^{2}y_2^{4}  \\
&&  +y_1^{4}y_2^{4}-y_1^{5}y_2^{4}-y_1^{4}y_2^{5}
+y_1^{6}y_2^{6},
\end{eqnarray*}
and the denominator is
\begin{eqnarray*}
D(1,2) &=& (1-x_1)(1-x_1^2)(1-y_1^2)(1-y_1y_2)(1-y_2^2)\cdot \\
&& (1-x_1y_1^2)(1-x_1y_1y_2)(1-x_1y_2^2)(1-y_1^2y_2)
(1-y_1y_2^2)\cdot \\
&& (1-x_1^2y_1^2)(1-x_1^2y_1y_2)(1-x_1^2y_2^2)\cdot \\
&& (1-y_1^4)(1-y_1^3y_2)(1-y_1^2y_2^2)(1-y_1y_2^3)(1-y_2^4).
\end{eqnarray*}

The numerator $N(0,3)$ has 3188 terms and so we shall omit it.
However, if we set $y_1=y_2=y_3=t$ then we obtain the simply graded
Poincar\'{e} series $P(\Cb;0,3;t)$. After these substitutions
the numerator and the denominator acquire common factors. We cancel
these factors and obtain normalized numerator and denominator
which we denote simply by $N(0,3;t)$ and $D(0,3;t)$.
The numerator is the palindromic polynomial
\begin{eqnarray*}
N(0,3;t) &=& 1-t-{t}^{2}-3\,{t}^{3}+18\,{t}^{4}+2\,{t}^{5}+12\,{t}^{6}
-36\,{t}^{7}+62\,{t}^{8} \\
&& +35\,{t}^{9}+124\,{t}^{10}-60\,{t}^{11}+81\,{t}^{12}+15\,{t}^{13}
+234\,{t}^{14} \\
&& +15\,{t}^{15}+\cdots-{t}^{27}+{t}^{28},
\end{eqnarray*}
and the denominator
$$ D(0,3;t)=\left( 1-t \right)  \left( 1-{t}^{2} \right) ^{7}
 \left( 1-{t}^{3} \right) ^{5} \left( 1-{t}^{4} \right) ^{7}. $$
These formulae agree with the computations of B. Broer
\cite[Section 8.2]{BB} for the group $\SO_5$.

\section{The algebra $\Cb$ in the case $n=2$ and $k_1=k_2=2$}
\label{2-mat-4}

In this case $\Cb$ is the algebra of symplectic polynomial invariants of 
two generic complex $4\times4$ matrices which we denote by $x$ and $y$. 
Our objective in this section is to find an MSG for $\Cb$.
We were not able to compute the multigraded Poincar\'{e} series
$P(\Cb;x_1,x_2,y_1,y_2)$. Thus, in order to simplify the computation, we set
$x_1=y_1=s$ and $x_2=y_2=t$. We refer to this new rational function
as the \emph{bigraded Poincar\'{e} series} of $\Cb$ and denote it by
$P(\Cb;s,t)$. We use similar notation for its numerator and
denominator.

The numerator $N(\Cb;s,t)$ has 284 terms and total degree 46.
It suffices to record only the terms of total degree $\le23$
because it satisfies the identity
$$ N(\Cb;s,t)=s^{23}t^{23}N(\Cb;s^{-1},t^{-1}). $$
The numerator is given by

\newpage
\begin{eqnarray*}
&& N(\Cb;s,t)= \\
&& 1+2st(s^{2}+3st+t^{2})+st(s+t)(3s^{2}+5st+3t^{2}) \\
&& +st(2s^{4}+9s^{3}t+12s^{2}t^{2}+9st^{3}+2t^{4}) \\
&& +st(s+t)( s^{4}+3s^{3}t+6s^{2}t^{2}+3st^{3}+t^{4} ) \\
&& +s^{2}t^{2}( s^{4}+7s^{3}t+17s^{2}t^{2}+7st^{3}+t^{4} ) \\
&& +s^{3}t^{3}( s+t )( 9s^{2}+13st+9t^{2} ) \\
&& +s^{3}t^{3} ( 9s^{4}+31s^{3}t+41s^{2}t^{2}+31st^{3}+9t^{4} ) \\
&& +s^{2}t^{2}(s+t)^{3}(s^{4}+4s^{3}t+8s^{2}t^{2}+4st^{3}+t^{4}) \\
&& +s^{3}t^{3}(2s^{6}+s^{5}t+9s^{4}t^{2}+21s^{3}t^{3}+9s^{2}t^{4}+st^{5}+2t^{6}) \\
&& -s^{3}t^{3}(s+t)(2s^{6}+11s^{5}t+16s^{4}t^{2}+15s^{3}t^{3}
+16s^{2}t^{4}+11st^{5}+2t^{6}) \\
&& -s^{3}t^{3}(2s^{8}+19s^{7}t+45s^{6}t^{2}+68s^{5}t^{3}+74s^{4}t^{4}+\cdots ) \\
&& -s^{4}t^{4}(s+t)(9s^{6}+23s^{5}t+37s^{4}t^{2}+40s^{3}t^{3}
+37s^{2}t^{4}+23st^{5}+9t^{6}) \\
&& -s^{4}t^{4}(2s^{8}+9s^{7}t+38s^{6}t^{2}+54s^{5}t^{3}+57s^{4}t^{4}+\cdots ) \\
&& +s^{4}t^{4}(s+t)(s^{8}+s^{7}t-10s^{6}t^{2}-20s^{5}t^{3}-21s^{4}t^{4}-\cdots ) \\
&& +s^{5}t^{5}(2s^{8}-7s^{7}t-27s^{6}t^{2}-52s^{5}t^{3}-57s^{4}t^{4}-\cdots ) \\
&& -4s^{6}t^{6}(s+t)^{3}(s^{2}+t^{2})(2s^{2}+st+2t^{2}) \\
&& -s^{6}t^{6} ( 4s^{8}+11s^{7}t+33s^{6}t^{2}+38s^{5}t^{3}+35s^{4}t^{4}+\cdots ) \\
&& +s^{7}t^{7}(s+t)(4s^{6}+11s^{5}t+21s^{4}t^{2}+36s^{3}t^{3}
+21s^{2}t^{4}+11st^{5}+4t^{6}) \\
&& +s^{6}t^{6}(s^{10}+13s^{9}t+45s^{8}t^{2}+94s^{7}t^{3}+152s^{6}t^{4}
+183s^{5}t^{5}+\cdots ) \\
&& +2s^{7}t^{7}(s+t)(3s^{8}+13s^{7}t+30s^{6}t^{2}+52s^{5}t^{3}+56s^{4}t^{4}+\cdots ) \\
&& +\cdots \\
&& +s^{19}t^{19}(s+t)(3s^{2}+5st+3t^{2})+2s^{20}t^{20}(s^{2}+3st+t^{2})+s^{23}t^{23}.
\end{eqnarray*}
The missing terms inside the parentheses should be inserted by using the obvious
symmetry. E.g. in the line beginning with $+s^4t^4(s+t)$ the missing part is
$-20s^3t^5-10s^2t^6+st^7+t^8$. The denominator is
\begin{eqnarray*}
D(\Cb;s,t) &=& (1-s)(1-t)(1-s^2)^2(1-st)^2(1-t^2)^2 \cdot \\
&& (1-s^3)(1-s^2t)^3(1-st^2)^3(1-t^3) \cdot \\
&& (1-s^4)^2(1-s^3t)^2(1-s^2t^2)^2(1-st^3)^2(1-t^4)^2.
\end{eqnarray*}

After setting $s=t$ and simplifying the fraction, the numerator becomes
the palindromic polynomial
\begin{eqnarray*}
N(\Cb;t) &=& 1-2\,t+{t}^{2}+2\,{t}^{3}+10\,{t}^{4}-4\,{t}^{5}
+7\,{t}^{6}+4\,{t}^{7}+52\,{t}^{8}  \\
&& +24\,{t}^{9}+38\,{t}^{10}+83\,{t}^{12}+56\,{t}^{13}+83\,{t}^{14}+
\cdots -2\,{t}^{25}+{t}^{26}
\end{eqnarray*}
and the denominator
$$ D(\Cb;t)=(1-t)^4(1-t^2)^6(1-t^3)^6(1-t^4)^6. $$

\begin{theorem} \label{MinSkup}
Let $\Cb$ be the algebra of $\Sp_4$-invariants of two generic
matrices $x,y\in M_4$. The  traces of the following $136$ matrices
form an MSG of $\Cb$:
\begin{eqnarray*}
&& x,\, y; \\
&& x^2,\, xx^*,\, xy,\, xy^*,\, y^2,\, yy^*; \\
&& x^3,\, x^2y,\, x^2y^*,\, xx^*y,\, y^2x,\, y^2x^*,\, yy^*x,\, y^3; \\
&& x^4,\, x^3x^*,\, x^3y,\, x^3y^*,\, x^2x^*y,\, x^*x^2y^*,\, x^2y^2,\, 
(xy)^2,\, xyxy^*,\, yxyx^*, \\
&& \quad xyx^*y^*, x^2yy^*,\, y^2xx^*,\, xx^*yy^*,\, y^3x,\, y^3x^*,\, 
y^2y^*x,\, y^*y^2x^*, y^4,\, y^3y^*; \\
&& x^3x^*y,\, x^2x^*xy,\, xx^*x^2y,\, x^3y^2,\, x^3yy^*,\, x^3y^*y,\,
x^2x^*y^2,\, x^*x^2y^2, \\ 
&& \quad x^2yx^*y^*,\, x^2y^*xy,\, xx^*yx^*y,\, y^3x^2,\, y^3xx^*,\,y^3x^*x,\, y^2y^*x^2,\, 
y^*y^2x^2, \\ 
&& \quad y^2xy^*x^*,\, y^2x^*yx,\, yy^*xy^*x,\, y^3y^*x,\, y^2y^*yx,\, yy^*y^2x; \\
&& x^2x^*x^2y,\, x^2(x^*)^2xy,\, x^2(x^*)^2y^2,\, x^*x^2x^*y^2,\, (x^2y)^2,\, x^*x^2yxy^*, \\
&& \quad x^3yxy^*,\, x^2yx^2y^*,\, x^2yxx^*y,\, x^2yx^*xy,\, x^2yx^*xy^*,\, 
x^3y^2y^*,\, y^3x^2x^*, \\
&& \quad x^2y^2xy,\, y^2x^2yx,\, x^2y^2xy^*,\, y^2x^2yx^*,\, x^2y^*yxy,\, y^2x^*xyx, \\
&& \quad xx^*(y^*)^2xy,\, yy^*(x^*)^2yx,\, xyx^*yx^*y^*,\, xyx^*y^*xy^*,\, y^2(y^*)^2x^2,\\
&& \quad  y^*y^2y^*x^2,\, (y^2x)^2,\, y^*y^2xyx^*,\, y^3xyx^*,\, y^2xy^2x^*,\, y^2xyy^*x, \\
&& \quad y^2xy^*yx,\, y^2xy^*yx^*,\, y^2y^*y^2x,\, y^2(y^*)^2yx; \\
&& x^2x^*x^3y,\, x^3x^*xy^2,\, xx^*xyx^2y,\, x^3yx^*xy^*,\, x^2x^*yx^2y^*,\, x^3y^2xy, \\
&& \quad x^*xx^*y^2xy^*,\, x(x^*)^2y^2x^*y,\, x^2yy^*x^2y,\, xx^*y^2x^2y,\, x^*xy^2x^2y, \\
&& \quad x^*xy^2xx^*y,\, x^2yxyxy^*,\, x^2yxyx^*y,\, y^3x^2yx,\, y^*yy^*x^2yx^*, \\
&& \quad y(y^*)^2x^2y^*x,\, y^2xx^*y^2x,\, yy^*x^2y^2x,\, y^*yx^2y^2x,\, y^*yx^2yy^*x, \\
&& \quad y^2xyxyx^*,\, y^2xyxy^*x,\, y^3y^*yx^2,\, yy^*yxy^2x,\, 
y^3xy^*yx^*,\, y^2y^*xy^2x^*, \\
&& \quad y^2y^*y^3x; \\
&& x^3yx^2x^*y,\, x^3y^2x^2y,\, x^*x^2y^*yx^2y,\, xx^*yx^2yxy,\, xx^*xy^3xy,\, yy^*yx^3yx, \\
&& \quad x^3y^3x^*y,\, y^3x^3y^*x,\, y^3x^2y^2x,\, y^*y^2x^*xy^2x,\, 
yy^*xy^2xyx,\, y^3xy^2y^*x; \\
&& x^3yx^2yxy,\, x^2y^2xyx^2y,\, y^2x^2yxy^2x,\, y^3xy^2xyx. \\
\end{eqnarray*}
\end{theorem}
\noindent {\bf Proof.} 
The proof is of computational nature. We give some details and 
describe what kind of computations we have actually performed in 
order to prove the theorem.

Let us begin by displaying in Table 1 the coefficients of the 
Taylor expansion
$$ P(\Cb;s,t)=\sum_{i,j=0}^\infty c_{i,j} s^it^j. $$
We need these coefficients only for $i+j\le11$. 
Since $c_{i,j}=c_{j,i}$, it suffices to record them for $i\le j$.

\begin{center}

{\bf Table 1: The coefficients $c_{i,j};\ i+j=k\le11$}
\[
\begin{array}{l|rrrrrrrr} k \diagdown i 
& 0 & 1 & 2 & 3 & 4 & 5 \\ \hline
 0& 1 &&&&&\\
 1& 1 &&&&&\\
 2& 3& 3 &&&&\\
 3& 4& 8 &&&&\\
 4& 9& 17& 28 &&&\\
 5& 11& 33& 59 &&&\\
 6& 20& 58& 133& 156 &&\\
 7& 25& 97& 238& 359 &&\\
 8& 41& 153& 437& 730& 906 &\\
 9& 50& 233& 703& 1372& 1907 &\\
 10& 75& 342& 1143& 2398& 3806& 4335 \\
 11& 91& 489& 1707& 3978& 6867& 8942 
\end{array} \]
\end{center}

Denote by $W$ the set of 136 words listed in the theorem.
Let $\pA$ denote the bigraded subalgebra of $\Cb$ generated by the traces
of the words in $W$. Let $\Cb_{i,j}$ resp. $\pA_{i,j}$
denote the homogeneous component of bidegree $(i,j)$ of $\Cb$ resp. $\pA$.
We know that $\dim \Cb_{i,j}=c_{i,j}$ for all $i,j$ and, of course,
$\pA_{i,j}\subseteq\Cb_{i,j}$.

Note that $\pA_{0,0}=\Cb_{0,0}$ is
the 1-dimensional space spanned by the constant 1. We also have
$\pA_{1,0}=\Cb_{1,0}$, the 1-dimensional space spanned by $\tr x$.
Similarly, $\pA_{0,1}=\Cb_{0,1}$ is the 1-dimensional space spanned by
$\tr y$.
According to Table 1, the subspace $\Cb_{2,0}$ has dimension 3.
It is easy to verify that the polynomials $(\tr x)^2$, $\tr x^2$ and
$\tr xx^*$ are linearly independent, and so they form a basis of
$\Cb_{2,0}$. It follows that $\pA_{2,0}=\Cb_{2,0}$. Similarly,
$\Cb_{1,1}$ is spanned by $(\tr x)(\tr y)$, $\tr xy$ and $\tr xy^*$,
and $\Cb_{0,2}$ is spanned by $(\tr y)^2$, $\tr y^2$ and $\tr yy^*$.
Hence $\pA_{1,1}=\Cb_{1,1}$ and $\pA_{0,2}=\Cb_{0,2}$.

To continue, it is convenient to use a computer. We used MAPLE 
with the ``LinearAlgebra'' package and its subpackage called 
``Modular''. For each bidegree $(i,j)$ with $i+j\le10$ we have 
constructed a set of $c_{i,j}$ linearly independent products 
$(\tr w_1)\cdots(\tr w_k)$, $w_r\in W$, of bidegree $(i,j)$. 
This means that we have verified the equality $\pA_{i,j}=\Cb_{i,j}$ 
for all indices $i$ and $j$ such that $i+j\le10$.
Since Kuzmin's conjecture is valid for $m=4$, the algebra $\Cb$
is generated by its homogeneous components of degree at most 
$10$, and so $\pA=\Cb$. As an additional check of the correctness 
of our result we have verified that the equalities 
$\pA_{i,j}=\Cb_{i,j}$ also hold for $i+j=11$.

The computations were actually carried out by using the compact 
real form $\Sp(2)$ of $\Sp_4$ and its diagonal conjugation action 
on the pair of generic quaternionic $2\times2$ matrices 
$x=(x_{ij})$ and $y=(y_{ij})$.
By using Schur's theorem, we may assume that $x$ is upper
triangular, i.e., $x_{21}=0$. By simultaneously
conjugating $x$ and $y$ with a diagonal matrix in $\Sp(2)$
we may also assume that $y_{2,1}\in\bR$, $y_{12}\in\bC$ 
and $x_{11}$ is in the $\bR$-span of the quaternionic units $1,i,j$. 
These simplifications were used in order to economize the memory 
and speed up the computation.

In the cases $i+j=11$ we had to use yet another simplification. 
Namely, we assumed that $\tr x=\tr y=0$, in which case the 
Poincar\'{e} series becomes 
$$ (1-s)(1-t)P(\Cb,s,t)=\sum_{i,j=0}^\infty c_{i,j}^0 s^i t^j, $$
and the dimensions $c_{i,j}^0$ are considerably smaller 
than $c_{i,j}$.
$\blacksquare$

\begin{remark} 
The maximum total degree of our generators is 9 while
the corresponding maximum degree for an MSG of $\SL_4$-invariants 
of two complex matrices is $10$.
\end{remark}

\begin{remark} 
Since the coefficients in the expansion of 
$(1+t+t^2)^2 N(\Cb;t)$ are nonnegative and
$$ (1+t+t^2)^2 D(\Cb;t)=(1-t)^2(1-t^2)^6(1-t^3)^8(1-t^4)^6, $$
it is likely that $\Cb$ has an HSOP consisting of
the traces of the first 16 words listed in the theorem and 
additional 6 invariants of degree 4. However, since $N(\Cb;s,t)$
vanishes for $s=t=1$, $\Cb$ does not have a bigraded HSOP.
\end{remark}

Let $W$ be defined as in the above proof. Denote by $\Cb(k)$ 
the subalgebra of $\Cb$ generated by the traces of the words 
$w\in W$ of length $\le k$.
Denote by $\Cb_d$ resp. $\Cb(k)_d$ the homogeneous component
of degree $d$ of the algebra $\Cb$ resp. $\Cb(k)$ with respect
to the $\bZ$-gradation by the total degree. The standard action of
$\GL_2$ on the 2-dimensional space spanned by the two generic
matrices $x$ and $y$ extends to an action on $\Cb$ by algebra
automorphisms preserving the $\bZ$-gradation.
Each homogeneous component $\Cb(k)_d$ is a polynomial
$\GL_2$-module. Let us introduce the quotient modules
$$ V_d=\Cb_d/\Cb(d-1)_d, \quad d=1,2,\ldots $$
Clearly these modules are nonzero only for $1\le d\le9$.
Each of our generators $\tr w$, $w\in W$, is a weight vector for 
the action of $\GL_2$ and its maximal torus consisting of the
diagonal matrices $\diag(\lambda,\mu)$. Hence it is easy to
determine the modules $V_d$ up to isomorphism.

Let $S_{p,q}$ denote a simple $\GL_2$-module with highest weight
$\lambda^p\mu^q$, $p\ge q\ge0$. If $S_{p,q}$ occurs as a submodule
of $V_d$, then necessarily $p+q=d$ and $d-p$ is even. Hence
two simple submodules of $V_d$ of the same dimension are
isomorphic. We shall write $V_d\simeq 1^{m_1}2^{m_2}3^{m_3}\cdots$
to indicate that $V_d$ is a direct sum of $m_1$ copies of a
1-dimensional module, $m_2$ copies of a 2-dimensional module, etc.
By using these conventions we can describe the structure of
these modules.

\begin{proposition} \label{modul}
The modules $V_d$ decompose into direct sum of simple
$\GL_2$-modules as follows:
\[ \begin{array}{lll}
 V_1\simeq 2^1 & \quad V_2\simeq 3^2 & \quad V_3\simeq 2^2\cdot 4^1 \\
 V_4\simeq 1^4\cdot 3^2\cdot 5^2 & \quad V_5\simeq 2^5\cdot 4^3
	 & \quad V_6\simeq 1^3\cdot 3^7\cdot 5^2 \\
 V_7\simeq 2^5\cdot 4^3\cdot 6^1 
	& \quad V_8\simeq 1^1\cdot 3^2\cdot 5^1 
	& \quad V_9\simeq 4^1
\end{array} \]
\end{proposition}
\noindent {\bf Proof.} 
This follows by inspecting the weights of $V_d$.
$\blacksquare$

We leave the job of computing a complementary $\GL_2$-submodule
of $\Cb(d-1)_d$ in $\Cb_d$ to the interested reader.

\section{Some Poincar\'{e} series $P(\Cb;\bx,\by)$ for
$n=3$ } \label{n=3}

In this section we assume that $n=3$, i.e., we work with 
the group $\Sp_6$. We give some multigraded Poincar\'{e} series
$P(\Cb;k_1,k_2)$, or their simplified expressions, for the algebra
$\Cb$ of $\Sp_6$-invariants of the  direct sum of $k_1$ copies of 
$\gp$ and $k_2$ copies of $\gk$ when $k_1+k_2\le3$.

It is well known that, for any $n$, the algebra $\Cb$ is regular when
$k_1+k_2=1$. Their Poincar\'{e} series are:
\begin{eqnarray*}
 P(\Cb;1,0) &=& \frac{1}{\prod_{k=1}^n (1-x_1^k)}, \\
 P(\Cb;0,1) &=& \frac{1}{\prod_{k=1}^n (1-y_1^{2k})}.
\end{eqnarray*}

Now let $n=3$ and $k_1+k_2=2$. If $k_1=2$ then $\Cb$ is a polynomial 
algebra in 10 generators, see \cite[Table 4a, case 5]{GS}.
For $k_1=2,1$ the numerators and the denominators
of $P(\Cb;k_1,k_2)$ are given by

\begin{eqnarray*}
N(2,0) &=& 1, \\
D(2,0) &=& (1-x_1)(1-x_2) (1-x_1^2)(1-x_1x_2)(1-x_2^2)  \cdot \\
	&& (1-x_1^3)(1-x_1^2x_2)(1-x_1x_2^2)(1-x_2^3) (1-x_1^2x_2^2), \\
N(1,1) &=& 1-x_1y_1(1+y_1^2)+x_1^2y_1^4(2+y_1^2)
	+x_1^3y_1^3(1+y_1+y_1^2-y_1^4) \\
	&& -x_1^4y_1^6(1+y_1)+x_1^5y_1^6(-1+y_1+y_1^2+y_1^3) \\
	&& +x_1^6y_1^7(1+2y_1^2)-x_1^7y_1^{10}+x_1^8y_1^{13}, \\
D(1,1) &=& (1-x_1)(1-x_1^2)(1-x_1^3)(1-y_1^2)(1-y_1^4)(1-y_1^6)\cdot \\
	&& (1-x_1y_1)(1-x_1y_1^2)(1-x_1y_1^3)(1-x_1y_1^4)\cdot \\
	&& (1-x_1^2y_1^2)^2(1-x_1^2y_1^4)(1-x_1^3y_1^2)(1-x_1^4y_1^2).
\end{eqnarray*}
Note that by setting $y_1=0$ resp. $x_1=0$ in $P(\Cb;1,1)$ we obtain
the rational function $P(\Cb;1,0)$ resp. $P(\Cb;0,1)$.

We have also computed the multigraded Poincar\'{e} series 
$P(\Cb;0,2)$ for the two copies of the adjoint module of $\Sp_6$.
Its numerator $N(0,2)$ has 508 terms. Fortunately there is a way 
to write it in a more compact form. We have
$$ N(0,2)=\sum_{d=0}^{66} c_d(y_1,y_2), $$
where $c_d(y_1,y_2)$ is a homogeneous polynomial of degree $d$.
The $c_d$s are symmetric polynomials, i.e., they satisfy 
$c_d(y_2,y_1)=c_d(y_1,y_2)$. Moreover, they satisfy the 
equation
$$ c_{66-d}(y_1,y_2)=(y_1y_2)^{33}c_d(\frac1{y_1},\frac1{y_2}). $$
By using the abbreviations
$$ \al=y_1y_2,\quad \be=y_1+y_2,\quad \ga=y_1^2+y_1y_2+y_2^2, $$
we have
\begin{eqnarray*}
N(0,2) &=& (1+\al^{33})-(\al+\al^{32})-(\al+\al^{31})\be
-(\al+\al^{30})(y_1^2+y_2^2) \\
&& -(\al+\al^{29})(y_1^3+y_2^3)+3(\al^2+\al^{29})\ga \\
&& +(\al^2+\al^{28})\be(2y_1^2+y_1y_2+2y_2^2) \\
&& +(\al^2+\al^{27})(3y_1^4+3y_1^3y_2+5y_1^2y_2^2+\cdots) \\
&& +(\al^2+\al^{26})\be(y_1^4+y_2^4) \\
&& +(\al^2+\al^{25})\ga(y_1^4-3y_1^3y_2-3y_1y_2^3+y_2^4) \\
&& -(\al^3+\al^{25})\be(2y_1+y_2)(y_1+2y_2)(y_1^2+y_2^2) \\
&& -3(\al^3+\al^{24})(y_1^6+y_1^5y_2+3y_1^4y_2^2+y_1^3y_2^3+\cdots) \\
&& -(\al^3+\al^{23})\be(y_1^6+2y_1^5y_2-y_1^4y_2^2
+2y_1^3y_2^3+\cdots) \\
&& -(\al^3+\al^{22})(y_1^8-y_1^7y_2-2y_1^6y_2^2-15y_1^5y_2^3
-8y_1^4y_2^4-\cdots) \\
&& +(\al^4+\al^{22})\be\ga(y_1^2+y_2^2)(y_1^2+6y_1y_2+y_2^2) \\
&& +(\al^4+\al^{21})(2y_1^8+5y_1^7y_2+18y_1^6y_2^2+14y_1^5y_2^3
+31y_1^4y_2^4+\cdots) \\
&& +(\al^4+\al^{20})\be(y_1^8+4y_1^7y_2+4y_1^6y_2^2+5y_1^5y_2^3
+5y_1^4y_2^4+\cdots) \\
&& +(\al^6+\al^{21})(2y_1^6-12y_1^5y_2-5y_1^4y_2^2
-25y_1^3y_2^3-\cdots) \\
&& -(\al^6+\al^{20})\be(4y_1^6+10y_1^5y_2+15y_1^4y_2^2
+22y_1^3y_2^3+\cdots) \\
&& -(\al^5+\al^{18})(y_1^{10}+5y_1^9y_2+21y_1^8y_2^2+24y_1^7y_2^3
+52y_1^6y_2^4+39y_1^5y_2^5+\cdots) \\
&& -(\al^6+\al^{18})\be(3y_1^8+8y_1^7y_2+9y_1^6y_2^2+17y_1^5y_2^3
+10y_1^4y_2^4+\cdots) \\
&& -(\al^6+\al^{17})(y_1^{10}+6y_1^9y_2+2y_1^8y_2^2+13y_1^7y_2^3
-8y_1^6y_2^4+7y_1^5y_2^5-\cdots) \\
&& +(\al^7+\al^{17})\be(y_1^8+7y_1^7y_2+8y_1^6y_2^2+24y_1^5y_2^3
+18y_1^4y_2^4+\cdots) \\
&& +(\al^7+\al^{16})\ga(y_1^2+y_2^2)(y_1^6+10y_1^5y_2+5y_1^4y_2^2
+19y_1^3y_2^3+\cdots) \\
&& +(\al^7+\al^{15})\ga(y_1^2+5y_1y_2+y_2^2)(y_1^3+y_2^3)
(y_1^4+y_1^3y_2+3y_1^2y_2^2+\cdots) \\
&& +(\al^8+\al^{15})(3y_1^{10}+3y_1^9y_2+14y_1^8y_2^2+6y_1^7y_2^3
+15y_1^6y_2^4-2y_1^5y_2^5+\cdots) \\
&& -(\al^9+\al^{15})\be(3y_1^8+4y_1^7y_2+15y_1^6y_2^2+16y_1^5y_2^3
+27y_1^4y_2^4+\cdots) \\
&& -(\al^9+\al^{14})\ga(4y_1^8+7y_1^7y_2+23y_1^6y_2^2+19y_1^5y_2^3
+34y_1^4y_2^4+\cdots) \\
&& -(\al^9+\al^{13})\be(2y_1^{10}+7y_1^9y_2+19y_1^8y_2^2 \\
&& \qquad\qquad\qquad\qquad +29y_1^7y_2^3+39y_1^6y_2^4
+37y_1^5y_2^5+\cdots) \\
&& -(\al^9+\al^{12})\ga(y_1^{10}+y_1^9y_2+8y_1^8y_2^2+9y_1^7y_2^3
+16y_1^6y_2^4+6y_1^5y_2^5+\cdots) \\
&& +(\al^{11}+\al^{13})\be\ga(y_1^6+2y_1^5y_2+3y_1^4y_2^2
+5y_1^3y_2^3+\cdots) \\
&& +(\al^{10}+\al^{11})(y_1^{12}+5y_1^{11}y_2+17y_1^{10}y_2^2 \\
&& \qquad\qquad\qquad\qquad +36y_1^9y_2^3+60y_1^8y_2^4
+74y_1^7y_2^5+81y_1^6y_2^6+\cdots) \\
&& +2\al^{11}\be\ga(y_1^8+3y_1^7y_2+7y_1^6y_2^2+8y_1^5y_2^3
+8y_1^4y_2^4+\cdots), \\
\end{eqnarray*}
where in some symmetric factors we suppress the terms 
$y_1^k y_2^l$ with $k<l$. 
(They can be easily recovered by using the symmetry.)
Note that we have combined the terms 
$c_d(y_1,y_2)$ and $c_{66-d}(y_1,y_2)$ for $0\le d<33$. 
For instance, the term $(\al^2+\al^{28})\be(2y_1^2+y_1y_2+2y_2^2)$ 
is the sum of 
$c_7(y_1,y_2)=\al^2\be(2y_1^2+y_1y_2+2y_2^2)$ and 
$ c_{59}(y_1,y_2)=\al^{28}\be(2y_1^2+y_1y_2+2y_2^2)$.

The denominator has a simple expression:
\begin{eqnarray*}
D(0,2) &=& (1-y_1^2)(1-y_1y_2)^2(1-y_2^2)\cdot \\
&& (1-y_1^2y_2)(1-y_1y_2^2)\cdot \\
&& (1-y_1^4)(1-y_1^3y_2)^2(1-y_1^2y_2^2)^2(1-y_1y_2^3)^2(1-y_2^4)
\cdot \\
&& (1-y_1^4y_2)(1-y_1^3y_2^2)(1-y_1^2y_2^3)(1-y_1y_2^4)\cdot \\
&& \prod_{k=0}^6 (1-y_1^{6-k}y_2^k).
\end{eqnarray*}

For the simplified expression obtained by setting $y_1=y_2=t$, the 
numerator is the palindromic polynomial
\begin{eqnarray*}
N(\Cb;t) &=& 1-t-2{t}^{3}+{t}^{4}-2{t}^{5}+12{t}^{6}-5{t}^{7}
+17{t}^{8}-20{t}^{9}  \\
&& +13{t}^{10}-35{t}^{11}+44{t}^{12}-26{t}^{13}+90{t}^{14}
-31{t}^{15}+90{t}^{16} \\
&& -87{t}^{17}+82{t}^{18}-106{t}^{19}+132{t}^{20}
-51{t}^{21} +170{t}^{22} \\
&& -51{t}^{23}+\cdots -{t}^{43}+{t}^{44}
\end{eqnarray*}
and the denominator is
$$ D(\Cb;t) = (1-t)(1-t^2)^3(1-t^3)^2(1-t^4)^7(1-t^5)^3(1-t^6)^5. $$

Finally we have computed the multigraded Poincar\'{e} series $P(\Cb;3,0)$:
\begin{eqnarray*}
 N(3,0) &=& 1+x_1x_2x_3(x_1x_2+x_2x_3+x_3x_1) \\
	&& +x_1^2x_2^2x_3^2(1+x_1^2x_2^2+x_2^2x_3^2+x_3^2x_1^2) \\
	&& +x_1^3x_2^3x_3^3( x_1+x_2+x_3 +x_1x_2+x_2x_3+x_3x_1 ) \\
	&& -x_1^4x_2^4x_3^4( x_1+x_2+x_3 +x_1x_2+x_2x_3+x_3x_1
	+x_1^2+x_2^2+x_3^2 ) \\
	&& -x_1^6x_2^6x_3^6(1+ x_1+x_2+x_3 )-x_1^8x_2^8x_3^8,  \\
 D(3,0) &=& (1-x_1)(1-x_2)(1-x_3)(1-x_1^2)(1-x_2^2)(1-x_3^2) \cdot \\
  	&& (1-x_1x_2) (1-x_2x_3) (1-x_3x_1) (1-x_1^3) (1-x_2^3)
	 (1-x_3^3) \cdot  \\ 
	&& (1-x_1^2x_2) (1-x_2^2x_3) (1-x_3^2x_1) \cdot \\
	&& (1-x_1x_2^2) (1-x_2x_3^2) (1-x_3x_1^2)(1-x_1x_2x_3) \cdot \\
	&& (1-x_1^2x_2^2) (1-x_2^2x_3^2) (1-x_3^2x_1^2) \cdot \\
	&&  (1-x_1^2x_2x_3)(1-x_1x_2^2x_3)(1-x_1x_2x_3^2).
\end{eqnarray*}

\section{The algebra $\Cb$ in the case $n=3$ and $k_1=k_2=1$}
\label{1-mat-6}

The bigraded Poincar\'{e} series $P(\Cb;1,1)$ for this algebra $\Cb$, as
given in the previous section, can be modified by multiplying the
numerator $N(1,1)$ and the denominator $D(1,1)$ with
$(1+x_1y_1+x_1^2y_1^2)(1+x_1y_1^3)$.
Then we obtain the polynomials
\begin{eqnarray*}
N^*(1,1) &=& 1+x_1^2y_1^4
+x_1^3 y_1^4 ( 1+y_1+y_1^2+y_1^3+y_1^5 ) \\
&& +x_1^4 y_1^4 ( 1+y_1 )  ( 1+y_1^2+y_1^4 )
+x_1^5 y_1^5 ( 1+y_1^2+y_1^3+y_1^4+y_1^5 ) \\
&& +x_1^6 y_1^8 ( 1+y_1+y_1^2+y_1^3+y_1^5 )
+x_1^7 y_1^9 ( 1+y_1 )  ( 1+y_1^2+y_1^4 ) \\
&& +x_1^8 y_1^9 ( 1+y_1^2+y_1^3+y_1^4+y_1^5 )
+x_1^9y_1^{14}+x_1^{11}y_1^{18}
\end{eqnarray*}
and
\begin{eqnarray*}
D^*(1,1) &=& (1-x_1)(1-x_1^2)(1-y_1^2)(1-x_1^3)(1-x_1y_1^2) \cdot \\
	&& (1-x_1^2y_1^2)^2(1-y_1^4)(1-x_1^3y_1^2)(1-x_1y_1^4) \cdot \\
	&& (1-x_1^4y_1^2)(1-x_1^3y_1^3)(1-x_1^2y_1^4)(1-y_1^6) (1-x_1^2y_1^6),
\end{eqnarray*}
respectively. All coefficients of $N^*(1,1)$ are nonnegative
integers (in fact they are all 0 or 1)
and $D^*(1,1)$ is still a product of terms $1-\mu$
where $\mu$ is a monomial in $x_1$ and $y_1$.
This suggests that $\Cb$ may have a bigraded HSOP.
We shall resolve this question in the next section.

We proceed to the (total degree) $\bZ$-gradation. Thus we set $x_1=y_1=t$
in $N^*(1,1)$ and $D^*(1,1)$ and we obtain the polynomials
\begin{eqnarray*}
N(\Cb;t) &=& 1+{t}^{6}+{t}^{7}+2{t}^{8}+2{t}^{9}+3{t}^{10}+{t}^{11}
+3{t}^{12}+2{t}^{13}  \\
&& +2{t}^{14}+2{t}^{15}+2{t}^{16}+3{t}^{17}+{t}^{18}
+3{t}^{19}+2{t}^{20}+2{t}^{21} \\
&& +{t}^{22}+{t}^{23}+{t}^{29}
\end{eqnarray*}
and
$$ D(\Cb;t)=(1-t)(1-t^2)^2(1-t^3)^2(1-t^4)^3(1-t^5)^2
(1-t^6)^4(1-t^8). $$
The Taylor expansion of $P(\Cb;t)=N(\Cb;t)/D(\Cb;t)$ is
\begin{eqnarray*}
P(\Cb;t) &=& \sum_{k\ge0} c_k t^k \\
&=& 1+t+3\,{t}^{2}+5\,{t}^{3}+11\,{t}^{4}+17\,{t}^{5}+35\,{t}^{6}+52\,{t}
^{7}+95\,{t}^{8} \\
&& +145\,{t}^{9}+245\,{t}^{10}+366\,{t}^{11}+597\,{t}^{12}+876\,{t}^{13}
+1368\,{t}^{14} \\ 
&& +1991\,{t}^{15}+3010\,{t}^{16}+4313\,{t}^{17}+6367\,{t}^{18}
+8992\,{t}^{19} \\
&& +12973\,{t}^{20}+18103\,{t}^{21}+25606\,{t}^{22}+35270\,{t}^{23}
+49072\,{t}^{24} \\
&& +66771\,{t}^{25}+91472\,{t}^{26}+123087\,{t}^{27}+166301\,{t}^{28} \\
&& +221356\,{t}^{29}+295393\,{t}^{30}+\cdots
\end{eqnarray*}
The coefficient $c_k$ gives the dimension of the homogeneous component
of $\Cb$ of degree $k$.

We used the same method as in section \ref{2-mat-4} to construct a candidate
MSG of $\Cb$. Let $z\in M_6$ be a generic matrix.
We introduce the set $W_1$ of 28 words in $z$ and $z^*$:
\begin{eqnarray*}
W_1 &=& \{ z;\, z^2,\, zz^*;\, z^3,\, z^2z^*;\, z^4,\, z^3z^*,\, z^2(z^*)^2;\,
z^5,\, z^4z^*; \\
&&  z^6,\, z^5z^*,\, z^4(z^*)^2,\, (z^2z^*)^2,\, z^2(z^*)^2zz^*;\, z^5(z^*)^2;\, \\
&& z^5z^*zz^*,\, z^4z^*z^2z^*,\, z^4z^*z(z^*)^2;\, z^5z^*z^2z^*,\, z^5z^*z(z^*)^2;\, \\
&& z^5z^*z^3z^*,\, z^5z^*z^2(z^*)^2,\, z^4z^*z^3(z^*)^2;\, z^4z^*z^3z^*zz^*;\, \\ 
&& z^5z^*z^4(z^*)^2,\, z^5z^*z^3z^*zz^*;\, z^5z^*z^4z^*zz^* \}
\end{eqnarray*}
and the set of their traces
$$ \tr W_1 = \{\tr w:\, w\in W_1\}. $$

We shall see below that $\tr W_1$ is an MSG of $\Cb$. One way to prove this
would use a bound for the maximum degree of generators in an MSG of $\Cb$.
However, the known bounds are too large.
The maximum degree of the polynomials in $\tr W_1$ is 13 and we
have verified that there are no additional generators in degrees $\le19$.
By replacing the generic matrix $z$ by one of trace 0 the computations
simplify somewhat and we were able to verify that there are no 
additional generators in degrees $20,\ldots,24$ as well, 
exceeding the (hypothetical) bound 21 derived from Kuzmin's conjecture.
The computations were performed by using the compact real form
$\Sp(3)$ of $\Sp_6$ and a generic matrix $z\in M_3(\bH)$. 
We can prove our assertion by using the $\bZ^2$-gradation of $\Cb$, avoiding
the use of Kuzmin's conjecture.
 
\begin{theorem} \label{ZgradMSG}
The set $\tr W_1$ is an MSG of $\Cb$.
\end{theorem}
\noindent {\bf Proof.} 
Recall that any two MSG's must contain the same number of polynomials
of degree $d$ for each $d$, see e.g. \cite[Algorithm 2.6.1]{DK}. 
The assertion of the theorem follows from 
Theorem \ref{msg} which we prove in the next section
and the computations mentioned above.
$\blacksquare$

Define the 15-element subset $W_1^\#\subseteq W_1$ by 
\begin{eqnarray*}
W_1^\# &=& \{ z;\, z^2,\, zz^*;\, z^3,\, z^2z^*;\, z^4,\, z^3z^*,\, 
z^2(z^*)^2;\, z^5,\, z^4z^*; \\
&& z^6,\, z^5z^*,\,  z^4(z^*)^2,\, (z^2z^*)^2;\, z^5z^*zz^* \}.
\end{eqnarray*}
Note that the degrees of the polynomials in $\tr W_1^\#$ match those
of the factors $1-t^d$ in the expression for the denominator $D(\Cb;t)$.
Let us denote by $i_k$ the generator $\tr w(z,z^*)$ of $\Cb$, where
$w(z,z^*)$ is the $k$th word in the listing defining $W_1$. 
For instance, we have $i_{10}=\tr z^4z^*$. By using this notation, we have
$$ \tr W_1^\# = \{ i_1,i_2,\ldots,i_{13},i_{14},i_{17} \}. $$

We have used in our computations the set $\tr W_1^\#$ as a tentative
HSOP of $\Cb$. If so, then $\Cb$ would be a free module
of rank 36 over the subalgebra generated by $\tr W_1^\#$. This module
would have a free basis consisting of homogeneous polynomials and
the coefficients of the numerator $N(\Cb;t)$ would give the number of
basis elements of each degree. We have constructed such a tentative
basis $I$, except for the last basis element of degree 29:
\begin{eqnarray*}
I &=& \{ 1,i_{15},i_{16},i_{18},i_{19},i_{20},i_{21},i_{22},i_{23},
i_{24},i_{25},i_{15}^2,i_{26}, \\
&& i_{27},i_{15}i_{16},i_{28},i_{15}i_{18},i_{15}i_{19},i_{15}i_{20},
i_{15}i_{21},i_{15}i_{23}, \\
&& i_{15}i_{24},i_{15}i_{25},i_{16}i_{22},i_{16}i_{23},i_{15}i_{26},
i_{15}i_{28},i_{16}i_{26},i_{16}i_{27}, \\
&& i_{18}i_{26},i_{19}i_{26},i_{18}i_{28},i_{19}i_{28},i_{22}i_{26},
i_{25}i_{26},\, ? \}.
\end{eqnarray*}

\section{A bigraded HSOP of $\Cb$}
\label{Hir}

In this section we continue the study of the algebra
$\Cb$ in the case $n=3$ and $k_1=k_2=1$ but here we exploit its
bigraded structure. As in the previous section, we work with the
real form of $\Cb$ obtained from the action of $\Sp(3)$ on 
$M_3(\bH)$. We write $z\in M_3(\bH)$ as $z=x+y$
where $x=x^*$ is hermitian and $y=-y^*$ a skew-hermitian matrix. 
Thus
\begin{equation} \label{ZbRaz}
x=\frac{1}{2}\left( z+z^* \right),\quad
y=\frac{1}{2}\left( z-z^* \right).
\end{equation}

Let $\Pi$ be the 15-dimensional real subspace of $M_3(\bH)$ 
consisting of matrices $z$ whose components $x$ and $y$ satisfy 
the following three conditions: (a) $x$ is a real diagonal matrix, 
(b) the diagonal entries of $y$ belong to $\bR i$ and (c) the off-diagonal entries of $y$ are in the $\bR$-span of the quaternions 
$1,i,j$. It is an easy exercise to show that every $\Sp(3)$-orbit 
in $M_3(\bH)$ has a representative $z\in\Pi$. 
Hence the restriction map from $\Cb$ to complex valued polynomial 
functions on $\Pi$ is injective. Therefore, for computational 
purposes, we may assume without any loss of generality that 
$z\in\Pi$.

The algebra $\Cb$ is generated by the traces of products 
$$ w(x,y)=x^{r_1}y^{s_1}\cdots x^{r_m}y^{s_m}. $$
The number of occurrencies of $x$ in this expression, $\sum r_i$, 
is the {\em $x$-degree} of $\tr w(x,y)$,  and $\sum s_i$ its 
{\em $y$-degree}. We refer to the sum of these two degrees as 
the {\em total degree}. 
The $x$- and $y$-degrees define the $\bZ^2$-gradation of $\Cb$, 
while the total degree defines the $\bZ$-gradation.
Since a homogeneous component for the $\bZ^2$-gradation has much
smaller dimension than the ambient $\bZ$-graded component, 
one can carry out the computations beyond the
limits reached in the previous section.

Define the set 
\begin{eqnarray*}
W_2 &=& \{ x;\, x^2,\, y^2;\, x^3,\, xy^2;\, x^2y^2,\, (xy)^2,\, y^4;\,
x^2yxy,\, xy^4; \\
&& (x^2y)^2,\, x^2y^2xy,\, xy^3xy,\, (xy^2)^2,\, y^6;\, x^2yxy^3;\, \\
&& x^2yxy^2xy,\, x^2y^4xy,\, (xy^3)^2;\, x^2y^3(xy)^2,\, x^2y^3xy^3;\, \\
&& x^2y^3xyx^2y,\, xy^3(xy)^3,\, xy^4xy^2xy;\, xy^3xy^2(xy)^2; \\
&& x^2y^4(xy)^3,\, xy^5xy^3xy;\, xy^5(xy)^2xy^2 \}
\end{eqnarray*}
consisting of 28 words in $x$ and $y$, and its 15-element subset
\begin{eqnarray*}
W_2^\# &=& \{ x;\, x^2,\, y^2;\, x^3,\, xy^2;\, x^2y^2,\, (xy)^2,\, y^4;\,
x^2yxy,\, xy^4; \\ 
&& (x^2y)^2,\, x^2y^2xy,\, xy^3xy,\, y^6;\, (xy^3)^2 \}.
\end{eqnarray*}
Note that the corresponding words in $W_1$ and $W_2$ have the same
length, and the same is true for $W_1^\#$ and $W_2^\#$. 

Let us introduce the notation for the generators of $\Cb$ by
writing $j_k$ for $\tr w(x,y)$ where $w(x,y)$ is the $k$th word
in the listing defining $W_2$. For instance, we have
$j_9=\tr x^2yxy$. By using this notation, we have
\begin{equation} \label{hsop-gen}
\tr W_2^\# = \{ j_1,j_2,\ldots,j_{13},j_{15},j_{19} \}.
\end{equation}
Denote by $\Cb^\#$ the subalgebra of $\Cb$ generated by $\tr W_2^\#$. 

Let $\vf:M_3(\bH)\to M_6$ be the embedding defined by
$$ \vf(z)=\left[ \begin{array}{cc} z_0 & -\bar{z}_1 \\ z_1 & \bar{z}_0 
\end{array} \right], $$
where $z=z_0+jz_1$ with $z_0,z_1\in M_3$. In this way we can view $M_6$
as the complexification of the real algebra $M_3(\bH)$.
We have $\tr z=\tr \vf(z)$,
the ordinary trace of the complex matrix $\vf(z)$. We also have
$\vf(z^*)=\vf(z)^*$, the symplectic adjoint of $\vf(z)$.
We shall write $z\in M_6$ as $z=x+y$, where $x$ and $y$ are defined by
(\ref{ZbRaz}) using the symplectic adjoint operation ${}^*$.

We claim that $\tr W_2^\#$ is an algebraically independent set. It
suffices to prove that the 15 functions $(\tr w)\circ\vf:M_3(\bH)\to\bR$,
$w\in W_2^\#$, are algebraically independent. This can be easily proved by
finding a point $z\in M_3(\bH)$ at which the Jacobian of this set
of functions has rank 15. For instance, we can take the point
$$ z=\left[ \begin{array}{ccc} -1+i&j&1+j\\j&1+i&i+j\\-1+j&i+j&i
\end{array} \right]. $$

The following fact, verified by extensive computer calculations, will
be used in the proof of the next theorem: The homogeneous components
of degree $d$ of the algebra $\Cb$ and the subalgebra
generated by $\tr W_2$ coincide for $d\le31$. In fact we have verified
that for each bidegree $(i,d-i)$, $0\le i\le d$ and $d\le31$,
the homogeneous components of these algebras coincide.

\begin{theorem} \label{hsop}
The set $\tr W_2^\#$ is an HSOP of the algebra $\Cb$.
\end{theorem}
\noindent {\bf Proof.} 
We have shown above that $\tr W_2^\#$ is an algebraically independent set.
Let $z\in M_6$ and let $z=x+y$ with $x^*=x$ and $y^*=-y$. Assume that
$\tr w(x,y)=0$ for $w\in W_2^\#$. To prove the theorem, we have to
show that the $\Sp_6$-orbit of $z$ is unstable (i.e., that 0 lies in
its closure). Equivalently, we have to show that $\tr w(x,y)=0$ for 
all nontrivial words $w$ (see \cite[Definition 2.4.6 and Lemma 2.4.5]{DK}). 
Note that we denote by ``$w$'' a word in two
noncommuting indeterminates. On the other hand ``$w(x,y)$'' denotes the
matrix obtained by evaluating $w$ at our pair of matrices $x$ and $y$.

Since $x^*=x$, $x$ is (symplectically) similar to $X\oplus X^T$ for
some $X\in M_3$, see e.g. \cite[Table III, line 4]{DPWZ}. 
As $x,x^2,x^3\in W_2^\#$, we have $\tr x=\tr x^2=\tr x^3=0$ and it
follows that $X^3=0$ and $x^3=0$. Since $x\in\gp$ and $\tr x=0$,
$x$ has the form
$$
x=\left[ \begin{array}{ccc|ccc} \xi_1&\xi_2&\xi_3&0&\eta_1&\eta_2\\
\xi_4&\xi_5&\xi_6&-\eta_1&0&\eta_3\\ 
\xi_7&\xi_8&-\xi_1-\xi_5&-\eta_2&-\eta_3&0\\ \hline
0&\zeta_1&\zeta_2&\xi_1&\xi_4&\xi_7\\ 
-\zeta_1&0&\zeta_3&\xi_2&\xi_5&\xi_8\\
-\zeta_2&-\zeta_3&0&\xi_3&\xi_6&-\xi_1-\xi_5 \end{array} \right].
$$

As $y^*=-y$, $\tr y^k=0$ for odd $k$. Since $y^2,y^4,y^6\in W_2^\#$,
we have $\tr y^2=\tr y^4=\tr y^6=0$ and we deduce that $y^6=0$. 
The simple Lie algebra $\sp_6$ has exactly 8
nilpotent adjoint orbits, see e.g. \cite[p. 82]{CM}. As their
representatives, we choose the following matrices:
\begin{eqnarray*}
&& 1) \left[ \begin{array}{ccc|ccc} 0&1&&&&\\&0&1&&&\\&&0&&&1\\ \hline
&&&0&&\\&&&-1&0&\\&&&&-1&0 \end{array} \right], \quad
2) \left[ \begin{array}{ccc|ccc} 0&1&&&&\\&0&&&1&\\&&0&&&1\\ \hline
&&&0&&\\&&&-1&0&\\&&&&&0 \end{array} \right], \\
&& 3) \left[ \begin{array}{ccc|ccc} 0&1&&&&\\&0&&&1&\\&&0&&&\\ \hline
&&&0&&\\&&&-1&0&\\&&&&&0 \end{array} \right], \quad
4) \left[ \begin{array}{ccc|ccc} 0&1&&&&\\&0&1&&&\\&&0&&&\\ \hline
&&&0&&\\&&&-1&0&\\&&&&-1&0 \end{array} \right], \\
&& 5) \left[ \begin{array}{ccc|ccc} 0&&&1&&\\&0&&&1&\\&&0&&&1\\ \hline
&&&0&&\\&&&&0&\\&&&&&0 \end{array} \right], \quad
6) \left[ \begin{array}{ccc|ccc} 0&&&1&&\\&0&&&1&\\&&0&&&\\ \hline
&&&0&&\\&&&&0&\\&&&&&0 \end{array} \right], \\
&& 7) \left[ \begin{array}{ccc|ccc} 0&&&1&&\\&0&&&&\\&&0&&&\\ \hline
&&&0&&\\&&&&0&\\&&&&&0 \end{array} \right], \quad
8) \quad 0.
\end{eqnarray*}
Clearly we may assume that $y$ is one of these matrices and so we have
eight cases to consider.

In each of these cases we shall prove that $\tr w(x,y)=0$ is valid
for all $w\in W_2$. Then 
by invoking the computational result mentioned above, it remains to show 
that $\tr w(x,y)=0$ for all words $w$ of length greater than $30$.
This will be accomplished by constructing enough nontrivial words $w$
such that $w(x,y)=0$.

Case 1). We claim that $\tr w(x,y)=0$ for all $w\in W_2$. 
This is proved by a routine computation for which we used MAPLE. From the 
equations $\tr xy^4=\tr (xy^3)^2=0$ we obtain that $\zeta_1=\zeta_2=0$. Then from 
equations $\tr xy^2=\tr xy^3xy=0$ we obtain that $\zeta_3=\xi_7=0$. Thus $x$ 
now depends on 10 variables only: the seven $\xi$'s and the three $\eta$'s.
Consider the polynomial ring (over the rationals) in these 10 variables and 
let $I$ be the ideal generated by the polynomials $\tr w(x,y)$, $w\in W_2^\#$.
We have computed a Groebner basis of $I$ and verified that all polynomials 
$\tr w(x,y)$, $w\in W_2$, belong to $I$. Consequently, our claim is proved.

By using the Groebner basis, it is easy to verify that
$\xi_4^2+\xi_8^2,\, \xi_8^3,\, \xi_5^6$ belong to $I$. Since all polynomials
in $I$ vanish on $x$, we conclude that $\xi_4=\xi_5=\xi_8=0$. This implies
that all products $x^ry^s$ with $r\in\{1,2\}$ and $s\in\{1,\ldots,5\}$
have the form
$$ \left[ \begin{array}{ccc|ccc} 0&*&*&*&*&*\\&0&&*&*&*\\&&0&*&*&*\\ \hline
&&&0&&\\&&&&0&\\&&&*&*&0 \end{array} \right], $$
where the blank entries are zeros. Consequently,
\begin{equation} \label{proiz}
x^{r_1}y^{s_1}x^{r_2}y^{s_2}x^{r_3}y^{s_3}x^{r_4}y^{s_4}=0
\end{equation}
for all $r_i\in\{1,2\}$ and $s_i\in\{1,\ldots,5\}$. Hence
$\tr w(x,y)=0$ for all words $w$ of length greater than $30$.

Case 2). We have $y^4=0$. From the equations
$\tr xy^2=\tr xy^3xy=0$ we obtain that $\zeta_1=\zeta_2=0$.
Next from the equations $\tr x^2y^2=\tr (xy)^2=\tr x^2y^2xy=0$ we obtain that
$\zeta_3=\xi_4=0$. Finally, the equation $\tr (x^2y)^2=0$ gives that
$\xi_6=0$ or $\xi_7=0$. In both cases one can verify that (\ref{proiz})
is valid for all $r_i\in\{1,2\}$ and $s_i\in\{1,2,3\}$. 

Cases 3) and 4). We have $y^4=0$ in the case 3) and $y^3=0$ in the case 4). 
From the equation $\tr xy^2=0$ we obtain that $\xi_7=0$. 
Then from the equations $\tr x^2y^2=\tr (xy)^2=0$ we obtain
that $\xi_8=\xi_4$. We now simplify the matrix $x$ by using these two facts.
We claim that $\tr w(x,y)=0$ for all $w\in W_2$. 
Consider the polynomial ring in 12 variables: the six $\xi$'s, 
three $\eta$'s and three $\zeta$'s. Denote by $I$ the ideal 
generated by the polynomials $\tr w(x,y)$ for $w\in W_2^\#$.
We computed a Groebner basis of $I$ (it has size 64) and we
verified that $\tr w(x,y)\in I$ for all $w\in W_2$, which proves 
our claim. 

Let $P$ be the set of words
\begin{eqnarray*}
&& y^3x; \\
&& y^2xyxy,\, y^2xyx^2yxy,\, y^2xyx^2yx^2y,\, y^2xyx^2y^2,\, y^2xy^2; \\
&& y^2x^2yxy,\, y^2x^2yx^2yxy,\, y^2x^2yx^2yx^2,\, y^2x^2yx^2y^2,\, y^2x^2y^2; \\
&& yxyxyxy,\, yxyxyx^2yxy,\, yxyxyx^2yx^2y; \\
&& yxyx^2yxy,\, yxyx^2yx^2yxy,\, yxyx^2yx^2yx^2; \\
&& yx^2yxyxyxy,\, yx^2yxyxyx^2y,\, yx^2yxyx^2; \\
&& yx^2yx^2yxyxy,\, yx^2yx^2yxyx^2,\, yx^2yx^2yx^2.
\end{eqnarray*}
For each $w\in P$ we have computed the matrix $w(x,y)$, reduced each of its
entries modulo $I$, and verified that the result was a zero matrix.
Thus we have $w(x,y)=0$ for all $w\in P$.
On the other hand let 
$$ w=y^{s_1}x^{r_1}\cdots y^{s_m}x^{r_m}, $$
with $r_i\in\{1,2\}$ and $s_i\in\{1,2,3\}$ in case 3) and $s_i\in\{1,2\}$ in case 4),
be a word whose length exceeds 30. Since the trace is
invariant under cyclic permutation of the factors, without any loss
of generality we may assume that $s_1\ge s_i$ for $1<i\le m$. It is not hard to 
show that such $w$ must have one of the words in $P$ as the initial subword.
Consequently $w(x,y)=0$, and so $\tr w(x,y)=0$ for all words $w$ whose length
exceeds 30.

Case 5). We have $y^2=0$. From $\tr (xy)^2=0$ we obtain that 
$\zeta_1^2+\zeta_2^2+\zeta_3^2=0.$
If $\zeta_1=\zeta_2=\zeta_3=0,$ then $(xy)^2=0$ and $xyx^2y=0$.
Consequently, $\tr w(x,y)=0$ for all nontrivial words $w$.
Assume now that at least one  $\zeta_i$ is not zero. 
If $A\in\SO_3$, the complex special orthogonal group, 
then $a=A\oplus A^T\in\Sp_6$ commutes with $y$. By replacing $x$
with $axa^{-1}$, where $A$ is suitably chosen, we may assume that
$\zeta_1=1$, $\zeta_2=0$, $\zeta_3=i$. From $\tr xyx^2y=0$ we obtain that 
$\xi_7=-2i\xi_1-\xi_3-i\xi_5.$ Then from $\tr x^2=0$ we obtain that
$$ \eta_1=\xi_1^2-\xi_3^2+\xi_5^2+\xi_1\xi_5+\xi_2\xi_4+\xi_6\xi_8
-2i\xi_1\xi_3-i\xi_3\xi_5-i\eta_3. $$
Next, from $\tr (x^2y)^2=0$ we obtain that $\xi_4=i\xi_6$.
Finally, from $\tr x^3=0$ we obtain that $\xi_3=-i\xi_1$. 
One can now verify that $x^iyx^jyx^ky=0$ for $i,j,k\in\{1,2\}$.
It follows that $\tr w(x,y)=0$ for all nontrivial words $w$.

Case 6). We have $y^2=0$. From $\tr (xy)^2=0$ we obtain that $\zeta_1=0$
and then deduce that $(xy)^2=0$. From $\tr (x^2y)^2=0$ we obtain that 
$\xi_8\zeta_2=\xi_7\zeta_3$ and deduce that $xyx^2y=0$. 
It follows that $\tr w(x,y)=0$ for all nontrivial words $w$.

Case 7). We have $y^2=0$ and it is easy to verify that also
$x^iyx^jy=0$ for $i,j\in\{1,2\}$. 
It follows that $\tr w(x,y)=0$ for all nontrivial words $w$.

Case 8). Since $x^3=0$ and $y=0$, 
it follows that $\tr w(x,y)=0$ for all nontrivial words $w$.
$\blacksquare$

Let $J$ be the following set of 36 bihomogeneous invariants:
\begin{eqnarray*}
J &=& \{ 1,j_{14},j_{16},j_{17},j_{18},j_{20},j_{21},j_{22},j_{23},
j_{24},j_{25},j_{26},j_{14}^2, \\
&& j_{27},j_{14}j_{16},j_{28},j_{14}j_{17},j_{14}j_{18},j_{14}j_{20},
j_{14}j_{21},j_{14}j_{22}, \\
&& j_{14}j_{23},j_{17}j_{20},j_{16}j_{23},j_{16}j_{24},j_{14}j_{26},
j_{16}j_{26},j_{18}j_{25},j_{14}j_{28}, \\
&& j_{18}j_{26},j_{17}j_{27},j_{17}j_{28},j_{18}j_{28},j_{22}j_{27},
j_{25}j_{26},j_{14}j_{25}j_{26} \}.
\end{eqnarray*}
The bidegrees of these invariants are pairwise distinct and coincide with
the bidegrees of the monomials which occur in the numerator $N^*(1,1)$.

\begin{theorem} \label{msg}
The above set $J$ is a free basis of the $\Cb^\#$-module $\Cb$ and the
set $\tr W_2$ is an MSG of $\Cb$. 
\end{theorem}
\noindent {\bf Proof.} 
It follows from Theorem \ref{hsop} that $\Cb^\#$ is a polynomial algebra
in the 15 variables (\ref{hsop-gen}). Also it follows that, as a 
$\Cb^\#$-module, $\Cb$ is a free module. The coefficient of $t^d$ in
the numerator $N(\Cb;t)$ gives the number of basis elements of degree d
of this free module. Let $L$ be the $\bZ$-graded $\Cb^\#$-submodule of $\Cb$
generated by $J$. We have verified that the homogeneous components
of $L$ and $\Cb$ of degree $d$ coincide for all $d\le30$. 
Hence $L=\Cb$ and the first assertion is proved.

Since the subalgebra of $\Cb$ generated by $\tr W_2$ contains $\Cb^\#$ and
$J$, it also contains $L$. Hence $\tr W_2$ generates $\Cb$. The minimality
of $\tr W_2$ is obvious from the procedure that we used to construct it,
see \cite[Algorithm 2.6.1]{DK}.
$\blacksquare$

The monomials in the generators $j_k$ of total degree $d$ are linearly
independent for $d<14$. This follows from the fact that the number of such 
monomials is equal to the coefficient $c_d$ of $t^d$ in the Taylor 
expansion of $P(\Cb,t)$ given in the previous section. On the other hand,
the number of such monomials of degree 14 is 1369 while $c_{14}=1368$.
Hence there exists a single nontrivial linear relation
among these monomials. It occurs in bidegree $(6,8)$. For quaternionic $3\times 3$ 
matrices $z=x+y$ of trace 0 this relation is:
\begin{eqnarray*}
&& 72j_{10}j_{8}j_{4}j_{2}-2304j_{12}^{2}j_{3}+288j_{10}j_{6}j_{5}j_{2}
+288j_{10}j_{9}j_{3}j_{2}-72j_{10}j_{5}j_{3}j_{2}^{2} \\
&& +192j_{10}j_{6}j_{4}j_{3}+192j_{13}j_{5}j_{4}j_{3}-24j_{7}j_{5}j_{4}j_{3}^{2}+
36j_{7}j_{5}^{2}j_{3}j_{2}+48j_{8}j_{6}j_{5}j_{4} \\
&& +576j_{16}j_{5}j_{3}j_{2}-144j_{6}j_{5}^{2}j_{3}j_{2}
-24j_{6}j_{5}j_{4}j_{3}^{2}-96j_{9}j_{8}j_{4}j_{3}-72j_{9}j_{5}j_{3}^{2}j_{2} \\
&& -288j_{9}j_{7}j_{5}j_{3}+864j_{9}j_{6}j_{5}j_{3}-144j_{9}j_{8}j_{5}j_{2}
-144j_{14}j_{5}j_{4}j_{3} \\
&& -72j_{14}j_{7}j_{3}j_{2}+144j_{14}j_{6}j_{3}j_{2}-32j_{15}j_{5}j_{4}j_{2}
+2j_{5}j_{4}j_{3}^{3}j_{2}+144j_{14}^{2}j_{2} \\
&& +2304j_{23}j_{7}-2304j_{23}j_{6}-32j_{15}j_{4}^{2}j_{3}-24j_{15}j_{7}j_{2}^{2}+
192j_{15}j_{9}j_{4} \\
&& +48j_{15}j_{6}j_{2}^{2}-96j_{15}j_{11}j_{2}+9j_{5}^{2}j_{3}^{2}j_{2}^{2}
+24j_{5}^{3}j_{4}j_{3}+72j_{10}^{2}j_{2}^{2}  \\
&& +2304j_{21}j_{9}-1152j_{13}j_{7}^{2}-48j_{10}j_{5}^{2}j_{4}-576j_{10}j_{9}j_{6}
-384j_{21}j_{4}j_{3} \\
&& -576j_{21}j_{5}j_{2}+288j_{13}j_{5}^{2}j_{2}-1152j_{13}j_{9}j_{5}
-384j_{13}j_{10}j_{4}+1152j_{13}j_{7}j_{6} \\
&& -3j_{7}j_{3}^{3}j_{2}^{2}-144j_{7}j_{6}j_{5}^{2}+576j_{11}j_{10}j_{5}
-12j_{11}j_{3}^{3}j_{2}-144j_{11}j_{5}^{2}j_{3} \\
&& +18j_{8}j_{5}^{2}j_{2}^{2}+16j_{8}j_{4}^{2}j_{3}^{2}-1152j_{16}j_{10}j_{2}
-384j_{16}j_{8}j_{4}-2304j_{16}j_{6}j_{5} \\
&& +192j_{16}j_{4}j_{3}^{2}+192j_{19}j_{5}j_{4}+288j_{19}j_{6}j_{2}
+6j_{6}j_{3}^{3}j_{2}^{2}-1152j_{17}j_{6}j_{3} \\
&& -144j_{17}j_{8}j_{2}+576j_{17}j_{7}j_{3}+36j_{14}j_{8}j_{2}^{2}
-18j_{14}j_{3}^{2}j_{2}^{2}+288j_{14}j_{7}j_{6} \\
&& +96j_{14}j_{10}j_{4}-576j_{14}j_{13}j_{2}+288j_{14}j_{11}j_{3}
+8j_{8}^{2}j_{4}^{2}+4608j_{16}^{2} \\
&& -1152j_{19}j_{11}+432j_{6}^{2}j_{5}^{2}+144j_{9}^{2}j_{3}^{2}+288j_{9}^{2}j_{8}
-288j_{17}j_{5}^{2}-576j_{17}j_{14}  \\
&& +2304j_{17}j_{13}-288j_{14}j_{6}^{2}+72j_{11}j_{8}j_{3}j_{2}
-36j_{8}j_{6}j_{3}j_{2}^{2}+18j_{8}j_{7}j_{3}j_{2}^{2} \\
&& +48j_{8}j_{7}j_{5}j_{4}+12j_{8}j_{5}j_{4}j_{3}j_{2}
+4608j_{18}j_{12}-2j_{4}^{2}j_{3}^{4}-36j_{10}j_{4}j_{3}^{2}j_{2} \\
&& +72j_{17}j_{3}^{2}j_{2}-72j_{14}j_{5}^{2}j_{2}
-2304j_{16}j_{9}j_{3}+576j_{14}j_{9}j_{5}=0.
\end{eqnarray*}
The components $x$ and $y$ are given by (\ref{ZbRaz}).
To obtain such relation valid for arbitrary $z\in M_3(\bH)$,
we can simply apply the above identity to the traceless matrix
$$ z-\frac{1}{6}(\tr z)I_3 = \left[ x-\frac{1}{6}(\tr x)I_3 \right] +y. $$
It is clear that the above identity is also valid
for matrices $z\in M_6$ of trace 0 provided we define $x$ and $y$ by
(\ref{ZbRaz}) using the symplectic adjoint operation ${}^*$.

\section{Conclusion}

The complex symplectic group $K=\Sp_{2n}$ acts by conjugation 
$(a,x)\to axa^{-1}$, $a\in K$, $x\in M_{2n}$, on the space $M_{2n}$ 
of complex matrices of order $2n$. This space decomposes as a 
direct sum of two $K$-submodules $M_{2n}=\gk\oplus\gp$, where 
$\gk$ is the adjoint module, i.e., the Lie algebra of $K$, and 
$\gp$ is the unique complementary submodule. For more details 
see Eqs. (\ref{def-k}) and (\ref{def-p}). We denote by $V$ 
the direct sum of $k_1$ copies of $\gp$ and $k_2$ copies of $\gk$. 
Let $\pP$ be the algebra of complex holomorphic polynomial 
functions on $V$, and $\Cb$ its subalgebra of $K$-invariant 
polynomials. One can use the above direct decomposition of $V$ 
to define a $\bZ^{k_1+k_2}$-gradation (multigradation) 
on $\pP$ and its subalgebra $\Cb$.
When $k_1=k_2=k$ we refer to $\Cb$ as the algebra 
of symplectic invariants of $k$ matrices 
$x_1,\ldots,x_k\in M_{2n}$. In that case the action of $K$ is 
given by Eq. (\ref{dij-akcija}).

It is well-known that, in the case of $k$ matrices, $\Cb$ is 
generated by the traces of all words in $x_1,\ldots,x_k$ and
their symplectic adjoints $x_1^*,\ldots,x_k^*$ (see the 
section \ref{Objasnjenja} for the definition).

First, let $\Cb$ be the algebra of symplectic invariants of two
matrices of order four (the case $n=k_1=k_2=2$). 
We have constructed an MSG 
(a minimal set of homogeneous generators) of $\Cb$. 
It consists of 136 traces (see Theorem \ref{MinSkup}). 
The degrees of these generators are $\le 9$.

Next, let $\Cb$ be the algebra of symplectic invariants of a 
single matrix of order six (the case $n=3$, $k_1=k_2=1$). 
For this algebra we have constructed an MSG consisting of 
28 traces (see Theorem \ref{ZgradMSG}). Moreover, we have 
constructed a polynomial subalgebra $\Cb^\#$ of $\Cb$, which is 
generated by 15 algebraically independent traces such that 
$\Cb$ is a free $\Cb^\#$-module of rank 36. 
(This is an example of Hironaka's decomposition.)

These results are based on extensive computations perfomed by
the MAPLE package for symbolic computations. As a preliminary 
but essential step, we have computed the Poincar\'{e} series 
for $\Cb$ as explicit rational expressions. In fact we have 
calculated these Poincar\'{e} series for the multigraded 
algebra $\Cb$ of the $K$-module $V$ 
when $n=2,3$ and for several choices of $k_1$ and $k_2$. 
The use of multigradation was necessary in a few cases in 
order to be able to carry out the required computations.
Indeed, the dimensions of the multigraded homogeneous 
components are much smaller than those of the ambient 
$\bZ$-graded homogeneous components. 

However we have not constructed a Hironaka's decomposition for 
the algebra of symplectic invariants of two matrices of order four. 
We leave this as an open problem.

\end{document}